\documentclass{amsart}
\usepackage{amsthm, amscd, amssymb, amsfonts,dsfont}
\usepackage{mathrsfs}
\usepackage[latin1]{inputenc}
\usepackage{graphicx}
\usepackage{enumerate}
\usepackage{esint}
\usepackage{color}
\usepackage{hyperref}

\newtheorem{thm}{Theorem}
\newtheorem{cor}{Corollary}
\newtheorem{lem}{Lemma}
\newtheorem{prop}{Proposition}
\newtheorem{rem}{Remark}

\newcommand{\N}{\mathds N}
\newcommand{\R}{\mathds R}
\newcommand{\fa}{\mathscr{F}}
\newcommand{\te}{\theta}
\newcommand{\p}{p(\cdot)}
\newcommand{\pri}{p'(\cdot)}
\newcommand{\be}{\beta}
\newcommand{\loc}{\text{\upshape loc}}
\newcommand{\Bro}{\mathcal{B}_\rho}

\begin{document}
	
\title{Weighted inequalities for Schr\"odinger type Singular Integrals on variable Lebesgue spaces}

\author{Adri\'an Cabral}


	
\thanks{This research is supported by  Universidad Nacional del Nordeste (UNNE) and Consejo Nacional de Investigaciones Cient\'ificas y T\'ecnicas (CONICET), Argentina.}

	
\address{Departamento de Matem\'atica -- Facultad de Ciencias Exactas y Naturales y Agrimensura -- UNNE and Instituto de Modelado e Innovaci\'on Tecnol\'ogica, CONICET--UNNE, Corrientes, Argentina.\\
Tel.: +54 (0379) 4471752}

\email{enrique.cabral@comunidad.unne.edu.ar}


\begin{abstract}
In this paper we study the boundedness in weighted variable Lebesgue spaces of  operators associated with the semigroup generated by the time-independent Schr\"odinger operator $\mathcal{L}=-\Delta+V$  in $\R^d$, where $d>2$ and the non-negative potential $V$ belongs to the reverse H\"older class $RH_q$ with $q>d/2$. Each of the operators that we are going to deal with are singular integrals  given by a kernel $K(x,y)$, which satisfies certain size and smoothness conditions in relation to a critical radius function $\rho$ which comes  appears naturally in the harmonic analysis related to  Schr\"odinger operator $\mathcal{L}$.

\end{abstract}

\maketitle
	
\section{Introduction and Preliminaries}
In this work we will consider operators  associated to the semi--group generated by the Schr\"odinger differential operator $\mathcal{L}=-\Delta+V$.

We will deal particularly with  first and second order Riesz transforms $R_1=\nabla\mathcal{L}^{-1/2}$ and  $R_2=\nabla^2\mathcal{L}^{-1}$ as  well as those involving the potential, namely, $V^{1/2}\mathcal{L}^{-1/2}$, $V\mathcal{L}^{-1}$ and $V^{1/2}\mathcal{L}^{-1}$.

The study of these operators under the assumptions that $V$ is a non-negative, non-identically zero and locally integrable function belonging to $RH_q$ for $q>d/2$ and $d>2$, was started by Shen in~\cite{Shen}, where he proves $L^p$ boundedness for all them.

We remind that $V\in RH_q$ means that there exists a constant $C$ such that
\begin{equation}\label{reverse}
\left(\frac{1}{|B|}\int_B V(y)^q\,dy\right)^{1/q} \leq \frac{C}{|B|}
\int_B V(y)\,dy,
\end{equation}
for every ball $B\subset \R^d$.

By H\"older inequality we can get that  $RH_q\subset RH_p$, for $q\ge p>1$. One remarkable feature about the $RH_q$ class is that, if $V\in RH_q$ for some $q>1$, then there exists $\epsilon>0$, which depends only on $d$ and the constant $C$ in~(\ref{reverse}), such that  $V\in RH_{q+\epsilon}$. Therefore, it is equivalent consider  $q>d/2$ or $q\ge d/2$.

After Shen's work several authors have dealt with boundedness results for such operators  acting on different spaces. For weighted $L^p$ spaces see for example~\cite{BHS-Classes},~\cite{BCH-Extrapol}, \cite{Tang},~\cite{Ly} and~\cite{BCH-Singular Int}, for the behaviour on suitable Hardy spaces
we refer to~\cite{DZ-H1},~\cite{Bui},~\cite{Liu-Tang-Zhu} and~~\cite{BCH-Singular Int}, as for the continuity on appropriate regularity spaces, see for instance~\cite{BHS-Riez Transform}, \cite{Ma-Stinga-et}, \cite{Bui} and~\cite{BCQ-Weighted STS}.

Under the above conditions on $V$, the analysis in the semi--group depends basically of a  function associated to $V$, given for $x\in\R^d$, by
\begin{equation}\label{rho Schrodinger}
\rho(x)=\sup\left\{r>0:\frac1{r^{d-2}}\int_{B(x,r)}V(y)\,dy\leq1\right\}.
\end{equation}
In particular, the reverse--H\"older condition implies that $0<\rho(x)<\infty$ for $x\in\R^d$.  Furthermore, according to \cite[Lemma 1.4]{Shen}, if $V\in RH_{q/2}$ the associated function $\rho$ verifies that there exist constants $c_\rho$, $N_\rho\geq1$ such that
\begin{equation}\label{est rho}
	c_\rho^{-1}\rho(x)\bigg(1+\frac{|x-y|}{\rho(x)}\bigg)^{-N_\rho}\leq\rho(y)\leq c_\rho\,\rho(x)\bigg(1+\frac{|x-y|}{\rho(x)}\bigg)^{\frac{N_\rho}{N_\rho+1}},
\end{equation}
for every $x,y\in\R^d$.

From now on, we call a \emph{critical radius function} to any positive  function $\rho$ that satisfies~(\ref{est rho}). Clearly, if $\rho$ is such a function, so it is $\beta\rho$ for any $\beta>0$.

Each of the operators that we are going to deal with are singular integrals  given by a kernel $K(x,y)$, which satisfies certain size and smoothness conditions in relation to the critical radius function $\rho$. These estimates of the kernel $K$ can be expressed in point-wise or integral way depending on the operator in question and the condition $RH_q$  that verifies the potential $V$.

So our strategy will be to work with general families of operators associated to a critical radius function $\rho$. When this function comes from a Schr\"odinger operator whose potential satisfies the aforementioned conditions, our families contains all the operators in which we are interested.

In this general context, we will prove our  boundedness results of such families of operators on weighted variable Lebesgue spaces $L^{\p}(w)$ with weights  defined in terms of $\rho$. Then, we apply those results to the Schr\"odinger setting to obtain bounding in $L^{\p}(w)$ for operators such that as first and second order Riesz transforms among others.

For this purpose, we developed extrapolation results that allow us to obtain boundedness results of the type described above in the variable setting. There is a natural bridge from the spaces $L^p(w)$ to $L^{\p}(w)$ given by the theory of Rubio de Francia extrapolation. Building on the results in~\cite{CU-Fiorenza-Neuge1}, Cruz-Uribe and Wang proved in~\cite{CUW} that if an operator $T$ maps $L^p(w)$ to itself whenever $w$ is in the Muckenhoupt $A_p$ class, then $T$ maps $L^{\p}(w)$ to itself all for all weights in $A_{\p}$ (see~\eqref{Ap punto} in Section~\ref{section:W and E}). In this paper we explore a similar connection between $L^p(w)$ and $L^{\p}(w)$ where the weights satisfy critical radius conditions $A_p^\rho$ and  $A_{\p}^\rho$ defined by~\eqref{A_p^rho} and~\eqref{Ap_rho} on Section~\ref{section:W and E}, but with a different approach. Instead we generalize an extrapolation theorem of the $L^\infty$--$BMO$ type proved by Harboure, Mac\'ias and Segovia in~\cite{Pola}. This requires to work with a weighted critical radius $BMO$ and critical radius sharp maximal operator, both introduced in~\cite{BCH-Extrapol}.

Our results cover a wide number of operators and, to our knowledge, the only result in this direction appeared in~\cite{Z-L}. In that article, the authors consider the Schr\"odinger type operator $L=-{\rm div}\,(A\nabla)+V$ on $\R^d$ with $d\ge3$, where the matrix $A$ satisfies uniformly elliptic condition  and prove that the operators $V L^{-1}$, $V^{1/2}\nabla L^{-1}$ and $\nabla^{2}L^{-1}$ are bounded on variable Lebesgue space $L^{\p}(\R^d)$ (without weights).

The remainder of this paper is organized as follows. In the rest of this section we present some basic results about variable Lebesgue spaces that will be used later. 

In the Section 2 we present the classes of weights involved in this work and some properties of them that will be useful. We also establish an important extrapolation theorem from the extreme $L^\infty$, which is a fundamental tool for the proof of our results.

In order to apply the extrapolation results of Section 2, we need to define a special sharp maximal function, since most of the interesting operators related to Schr\"odinger analysis are not bounded on $L^\infty$. We do this in Section 3 where we also present an appropriate version of $BMO$ space, introduced originally in~\cite{BCH-Extrapol}.

Section 4 is devoted to present two families of operators that include different singular integrals appearing in the Schr\"odinger setting and to establish and prove the main results of this paper. Finally we present some applications in the context of the Schr\"odinger operator in Section 5.

Throughout this paper, unless otherwise indicated, we will use $C$ and $c$ to denote
constants, which are not necessarily the same at each occurrence.

\

As we said before, we will now continue with some definitions related to the variable Lebesgue spaces. Let $\p:\R^d\rightarrow[1,\infty)$ be a measurable function. Given a measurable set $A\subset\R^d$ we define
\begin{equation*}
p^{-}_{A}=\text{ess}\inf_{x\in A}p(x),\hspace{1cm} p^{+}_{A}=\text{ess}\sup_{x\in A}p(x).
\end{equation*}

For simplicity we denote $p^{+}=p^{+}_{\R^d}$ and $p^{-}=p^{-}_{\R^d}$. Given $\p$, the conjugate exponent $\pri$ is defined pointwise
\begin{equation*}
\frac{1}{p(x)}+\frac{1}{p'(x)}=1,
\end{equation*}
where we let $p'(x)=\infty$ if $p(x)=1$.

By $\mathcal{P}(\R^d)$ we will denote  the collection of all measurable functions $\p:\R^d\rightarrow[1,\infty)$ and by $\mathcal{P}^*(\R^d)$ the set of $p\in\mathcal{P}(\R^d)$ such that $ p^{+}<\infty$.

Given  $p\in\mathcal{P}^*(\R^d)$, we say that a measurable function $f$ belongs to $L^{\p}(\R^d)$ if for some $\lambda>0$
\begin{equation*}
	\varrho(f/\lambda)=\int_{\R^d}\left(\frac{|f(x)|}{\lambda}\right)^{p(x)}dx<\infty.
\end{equation*}
A Luxemburg norm can be defined in $L^{\p}(\R^d)$ by taking
\begin{equation*}
	\|f\|_{\p}=\inf\{\lambda>0:\varrho(f/\lambda)\le1\}.
\end{equation*}

This spaces are special cases of Musielak-Orlicz spaces (see \cite{Musielak}), and generalize the classical Lebesgue spaces. For more information see,  for example, \cite{KR,Libro:CU-Fiorenza,Libro:D-H}.

We will denote with $L^{\p}_{\loc}(\R^d)$ the space of functions $f$ such that $f\chi_{B}\in L^{\p}(\R^d)$ for every ball $B\subset\R^d$.

The following conditions on the exponent arise in connection with the boundedness of the Hardy--Littlewood maximal operator $M$ in $L^{\p}(\R^d)$ (see, for example, \cite{Diening}, \cite{Libro:CU-Fiorenza} or~\cite{Libro:D-H}). We shall say that $p\in\mathcal{P}^{\log}(\R^d)$ if $p\in\mathcal{P}^*(\R^d)$ and if there are constants $C>0$ and $p_\infty$ such that
\begin{equation}
|p(x)-p(y)|\le \frac{C}{\log(e+1/|x-y|)},\hspace{1cm}\text{for every}~x,y\in\R^d,
\end{equation}
and
\begin{equation}
|p(x)-p_\infty|\le \frac{C}{\log(e+|x|)},\hspace{1cm}\text{with}~x\in\R^d.
\end{equation}
The functions in $\mathcal{P}^{\log}(\R^d)$ are generally called log--H\"older continuous.

Given a weight $w$ and $p\in\mathcal{P}(\R^d)$, following~\cite{CU-Fiorenza-Neuge}, define the weighted variable Lebesgue space $L^{\p}(w)$ to be the set of all measurable functions $f$ such that $fw\in L^{\p}(\R^d)$, and we write 
\begin{equation*}
\|f\|_{L^{\p}(w)}=\|f\|_{\p,w}=\|fw\|_{\p}.
\end{equation*} 
Thus, we say that an operator $T$ is bounded on $L^{\p}(w)$ if
\begin{equation*}
\|Tf\,w\|_{\p}\le C\|fw\|_{\p},
\end{equation*}
for all $f\in L^{\p}(w)$ and some constant $C>0$.

With these definitions  we have the following generalization of H\"older's inequality and an equivalent expression for the norm (see, for example, Section 2 in~\cite{Libro:CU-Fiorenza}).

\begin{lem}\label{lem: Holder}
Let $p,q,s\in\mathcal{P}(\R^d)$ be such that
\begin{equation*}
\frac1{s(x)}=\frac1{p(x)}+\frac1{q(x)}
\end{equation*}
for almost every $x$. Then, there exists a constant $C$ such that  
\begin{equation*}
\|fg\|_{s(\cdot)}\le C\|f\|_{p(\cdot)}\|g\|_{q(\cdot)},
\end{equation*}
for every $f\in L^{\p}(\R^d)$ and $g\in L^{q(\cdot)}(\R^d)$. In particular, $fg\in  L^{s(\cdot)}(\R^d)$.
\end{lem}

\begin{lem}\label{lem: norma dual}
Let $p\in\mathcal{P}(\R^d)$ and $w\in L^{\p}_{\loc}$. Then, there exist constants $c$ and $C$ such that
\begin{equation*}
c\|f\|_{L^{\p}(w)}\le\sup_{\|g\|_{L^{p'(\cdot)}(w^{-1})}\le1}\ \left|\int_{\R^d}f\,g\,dx\right|\le C\|f\|_{L^{\p}(w)}.
\end{equation*}
\end{lem}

\section{Weights and Extrapolation}\label{section:W and E}

In~\cite{BCH-Extrapol} the authors analyze the behavior of certain classes of operators that arise from the harmonic analysis related to the semigroup whose infinitesimal operator is the Schr\"odinger operator in weighted $L^p$ spaces for $1<p\le\infty$ for appropriate weights by means extrapolation techniques (see also~\cite{BCH-Singular Int}). To do this, it is essential to know the behavior of certain maximal operators in weighted Lebesgue spaces.

Given a critical radius function $\rho$, for each $\theta\ge0$ and $f\in L^1_{\loc}(\R^d)$ we define the maximal operator
\begin{equation}\label{Max_theta}
M^\te_{\rho} f(x) = \sup_{r>0}\bigg(1+\frac r{\rho(x)}\bigg)^{-\te} \frac{1}{|B(x,r)|} \int_{B(x,r)}|f(y)|\,dy\ \ \ \ \ x\in\R^d.
\end{equation}

By a weight we will mean a locally integrable function $w$ such that $0<w(x)<\infty$ almost everywhere.

We will also consider, as in~\cite{BHS-Classes}, a classes of weights associated to the critical radius function $\rho$, that properly contain the classical Muckenhoupt weights. For  $p>1$, the class $A_p^{\rho}$ is defined as those weights $w$ such that for some $\theta>0$, there exist a constant $C>0$ so that
\begin{equation}\label{A_p^rho}
\left(\int_{B}w\,dy\right)^{1/p}\left(\int_{B}w^{-\frac1{p-1}}\,dy\right)^{1/p'}\le C|B|\left(1+\frac{r}{\rho(x)}\right)^\theta,
\end{equation}
for every ball $B=B(x,r)$. Similarly, when $p=1$, we denote with $A_1^{\rho}$  the class of weights  $w$ such that for some $\theta>0$, there exist a constant $C>0$ so that
\begin{equation}\label{A_1^rho}
\int_{B}w\,dy\le C|B|\left(1+\frac{r}{\rho(x)}\right)^\theta\inf_{B}w,
\end{equation}
for every ball $B=B(x,r)$.

Classes $A_p^{\rho}$ are intimately connected with the family of maximal operators $M^\te_{\rho}$. In fact, for $1<p<\infty$ they are bounded on $L^p(w)$, provided $w\in A_p^{\rho}$ (see \cite[Proposition 3]{BCH-Extrapol}).

For our results we will need to know that the maximal operators $M^\te_{\rho}$ are bounded in  weighted variable Lebesgue spaces. The corresponding results about the boundedness of the the maximal Hardy--Littlewood function on  weighted  variable Lebesgue spaces  have been established in different settings by several authors, see for instance \cite{CU-Diening-Hasto,CU-Fiorenza-Neuge,Libro:D-H}. Following~\cite{CU-Fiorenza-Neuge},  if the weight $w$  is replaced  by $w^p$ in the definition of weights of Muckenhoupt, the condition $A_p$ for $p>1$ take the form
\begin{equation*}
\sup_{B}|B|^{-1}\|w\chi_B\|_{p}\|w^{-1}\chi_B\|_{p'}<\infty,
\end{equation*}
where the supreme is taken on all balls $B\subset\R^d$. 

This way of defining the weights $A_p$ immediately generalizes to variable Lebesgue spaces. Given an exponent $p\in\mathcal{P}(\R^d)$ and a weight $w$, it is said that $w\in A_{\p}$ if
\begin{equation}\label{Ap punto}
\sup_{B}|B|^{-1}\|w\chi_B\|_{\p}\|w^{-1}\chi_B\|_{p'(\cdot)}<\infty,
\end{equation}
where the supreme is taken on all balls $B\subset\R^d$.

Following the above and the ideas in~\cite{BCH-Extrapol} and~\cite{BHS-Classes},  in~\cite{C-Weighted norm} the authors defined classes of weights associated to a critical radius function $\rho$ in the variable context. Given a critical radius function $\rho$ and $p\in\mathcal{P}(\R^d)$ we will consider two families of weights. We introduce the  class of weights $A_{\p}^{\rho,\loc}$ as those $w$ satisfying~\eqref{Ap punto} for every ball $B\in\Bro$ with
\begin{equation}\label{def:Bro}
\mathcal{B}_\rho=\{B(x,r): x\in\R^d\ \text{and}\ r\le\rho(x)\}.
\end{equation}

Also,  we will say that a weight $w$ belongs to the class $A_{\p}^{\rho}$ if for some $\theta>0$, there exist a constant $C>0$ such that 
\begin{equation}\label{Ap_rho}
\|w\chi_B\|_{\p}\|w^{-1}\chi_B\|_{p'(\cdot)}\le C|B|\left(1+\frac{r}{\rho(x)}\right)^\te
\end{equation}
holds for all balls $B=B(x,r)$. Thus defined, it is clear that $A_{\p}\subset A_{\p}^\rho\subset A_{\p}^{\rho,\loc}$.

It is also clear that a weight $w$ belongs to the class $A_{\p}^{\rho}$ if and only if $w^{-1}$ belongs to the class $A_{\pri}^{\rho}$.

As shown in~\cite{C-Weighted norm} the weights $A_p^{\rho}$ and $ A_{\p}^{\rho,\loc}$ characterize the boundedness in $L^{\p}(w)$ of the  operators  $M^\te_{\rho}$ and  $M_{\rho}$, the latter defined for $f\in L^1_{\loc}(\R^d)$ as
\begin{equation*}\label{max_local}
M_\rho f(x) = \sup_{B\in\Bro}\frac{1}{|B|} \int_{B}|f(y)|\,dy.
\end{equation*}

\begin{thm}[\cite{C-Weighted norm}, Theorems~4 and 5]\label{teo: acot M}
Let $p\in\mathcal{P}^{\log}(\R^d)$ with $p^{-}>1$. Then, 
\begin{enumerate}[\ a)]
\item $w\in A_{\p}^{\rho,\loc}$ if and only if  $M_\rho$ is bounded in $L^{\p}(w)$.
\item $w\in A_{\p}^\rho$ if and only if there exists  $\te>0$ such that $M^\te_{\rho}$ is bounded in $L^{\p}(w)$.
\end{enumerate}
\end{thm}

We are now going to state and prove some facts about weights in the classes defined above which will be useful in what follows and are of interest by themselves.

\begin{prop}\label{prop: Aploc gamma}
Let  $p\in\mathcal{P}(\R^d)$ and $\be>0$, then $A_{\p}^{\beta\rho,\loc}=A_{\p}^{\rho,\loc}$. 
\end{prop}
\begin{proof}
This result was prove for $\beta>1$ in Proposition 3 of~\cite{C-Weighted norm}. When $\beta<1$, since $\varrho=\beta\rho$ is also a critical radius function and $1/\beta>1$, considering the previous case it follows that
\begin{equation*}
A_{\p}^{\beta\rho,\loc}=A_{\p}^{\varrho,\loc}=A_{\p}^{\frac1\beta\varrho,\loc}=A_{\p}^{\rho,\loc}.
\end{equation*}

\end{proof}

In order to prove the following result we need the following lemma, which was proved in~\cite{Libro:CU-Fiorenza}.

\begin{lem}[\cite{Libro:CU-Fiorenza}, Proposition 2.18]\label{lem:dilatacion}
Let $\p\in\mathcal{P}^*(\R^d)$ and $s$ such that $1/p^{-}\le s<\infty$. Then,
\begin{equation*}
\||f|^s\|_{\p}=\|f\|^s_{s\p}.
\end{equation*}
\end{lem}

\begin{prop}\label{prop: wdelta}
Let $0<\delta<1$ and $p\in\mathcal{P}^*(\R^d)$, with $p^{-}>\delta$. If $w\in A_{\p}^{\rho,\loc}$, then $w^\delta\in A_{\frac{\p}{\delta}}^{\rho,\loc}$.
\end{prop}
\begin{proof}
If $r(\cdot)=\p/(\p-\delta)$, we have  that $1/r(\cdot)=(1-\delta)+\delta/\pri$. Then, by the generalized H\"older inequality (Lemma~\ref{lem: Holder}), the above lemma and the hypothesis on $w$ and $\delta$, for every ball $B\in\Bro$, we have
\begin{equation*}
\begin{split}
\|w^{-\delta}\chi_B\|_{\big(\frac{\p}{\delta}\big)'}&=\|w^{-\delta}\chi_B\|_{\frac{\p}{\p-\delta}}\lesssim\|w^{-\delta}\chi_B\|_{\frac{\pri}{\delta}}\|\chi_B\|_{\frac1{1-\delta}}=\|w^{-1}\chi_B\|^\delta_{\pri}|B|^{1-\delta}.
\end{split}
\end{equation*}
	
Thus, given that $w\in A_{\p}^{\rho,\loc}$,  it follows
\begin{equation*}
\begin{split}
\|w^{\delta}\chi_B\|_{\frac{\p}{\delta}}\|w^{-\delta}\chi_B\|_{\big(\frac{\p}{\delta}\big)'}
&\le\|w\chi_B\|^\delta_{\p}\|w^{-1}\chi_B\|^\delta_{\pri}|B|^{1-\delta}\\
&\le C|B|^\delta|B|^{1-\delta}=C|B|.
\end{split}
\end{equation*}
\end{proof}

In the rest of this section we are interested in establishing some extrapolation results in the variable exponent spaces context.

Before presenting these results, it is necessary to establish the following technical lemmas.

\begin{lem}\label{algoritmo}
Let $r\in\mathcal{P}(\R^d)$, $S$ a sublinear operator and $\mu$  a weight such that $S$ is bounded on $L^{r(\cdot)}(\mu)$. Suppose further that $S(h)\ge0$ for any $h\in L^{r(\cdot)}(\mu)$. Then, if $h$ be a non-negative function such that  $h\in L^{r(\cdot)}(\mu)$ there exists a function $H$ such that $(1)$ $h\le H$ a.e.; $(2)$ $\|H\|_{L^{r(\cdot)}(\mu)}\le2\|h\|_{L^{r(\cdot)}(\mu)}$ and $(3)$ $S(H)\le C\,H$ a.e., with $C=2\|S\|_{L^{r(\cdot)}(\mu)}$.
\end{lem}
\begin{proof}
Following the algorithm of Rubio de Francia, it is sufficient to consider $H$ defined by
\begin{equation*}
H(x)=\sum_{k=0}^{\infty}\frac{S^{k}h(x)}{2^{k}\|S\|^{k}_{L^{r(\cdot)}(\mu)}},
\end{equation*}
where, for $k\ge1$, $S^{k}$ denotes $k$ iterations of the operator $S$, and $S^{0}$ is the identity operator. The properties of $H$ follow immediately from the hypothesis about $S$.

\end{proof}

\begin{lem}\label{lem:corlema4HMS}
Let $p\in\mathcal{P}(\R^d)$ and $w$ a weight. Then, for any non-negative function $f$ in $L^{p(\cdot)}(w)$ there exists a non-negative function $F$  such that
\begin{equation*}
\|F\|_{L^{p(\cdot)}(w)}\leq2
\end{equation*}
and
\begin{equation*}
\|f\|_{L^{p(\cdot)}(w)} = \|f F^{-1}\|_{L^\infty}.
\end{equation*}
\end{lem}

\begin{proof}
Let $H\in L^{p(\cdot)}(w)$ such that $|H(x)|>0$ for all $x\in\R^d$ and  $\|H\|_{L^{p(\cdot)}(w)}\le1$. Suppose  $\|f\|_{L^{p(\cdot)}(w)}=1$.
	
We consider
\begin{center}
$~ F(x)=\left\{\begin{array}{ccc}
f(x)&\text{if}&f(x)\neq0;\\
H(x)&\text{if}&f(x)=0.
\end{array}\right.$
\end{center}
	
By the definition of $F$ its follows that  $\|fF^{-1}\|_{L^\infty}=1=\|f\|_{L^{p(\cdot)}(w)}$.
	
On the other hand, we have
\begin{equation*}
\|F\|_{L^{\p}(w)}\le \|f\|_{L^{\p}(w)}+\|H\|_{L^{\p}(w)}\le2.
\end{equation*}
	
Now, for general $f$ general, we take $f/\|f\|_{L^{p(\cdot)}(w)}$ and we deduce
\begin{equation*}
\|f F^{-1}\|_{L^\infty}=\|f\|_{L^{\p}(w)}.
\end{equation*}
	
\end{proof}

Following~\cite{CU-Martell-Perez}, the extrapolation results of this section shall be expressed in terms of function pairs $(f,g)$ belonging to $\mathscr{F}$, a family of pairs of functions measurable and non-negative.

Given a weight $v$ such that $v^{-1}\in A_1^\rho$, the expression
\begin{equation}\label{des_inf}
\|f v\|_{L^\infty} \leq C \|g v\|_{L^\infty},\hspace{.5cm}(f,g)\in \mathscr{F},
\end{equation}
should be understood in the sense that the inequality holds for every $(f,g)\in \mathscr{F}$ whenever the left hand side is finite, with a constant $C$ depending on $v$ only through 
$[v^{-1}]_{1,\te}$, the infimum of the constants in~\eqref{A_1^rho} for $w=v^{-1}$.

\begin{thm}\label{thm:ext inf}
Let $\fa=\{(f,g)\}$  a family of pairs of measurable and non-negative functions. Suppose
\begin{equation}\label{des inf}
\|f w\|_{L^\infty} \le C \|g w\|_{L^\infty},
\end{equation}
holds for every weight $w$ such that $w^{-1} \in A_1^\rho$ and every pairs $(f,g)\in\mathscr{F}$. Then,
\begin{equation}\label{des infpp}
\|f\|_{L^{\p}(w)}\le C \|g\|_{{L^{\p}(w)}},
\end{equation}
holds for every pair $(f,g)\in\fa$, every $p\in\mathcal{P}^{\log}(\R^d)$ with $p^{-}>1$ and every $w\in A_{\p}^{\rho}$, provided that the left hand side of \eqref{des infpp} is finite.
\end{thm}	
\begin{proof}
Let $w\in A_{\p}^{\rho}$ and $(f,g)\in\mathscr{F}$. We may suppose, without loss of generality, that $f$ and $g$ belong to $L^{\p}(w)$. From Lemma~\ref{lem:corlema4HMS}, there exist non-negative functions $F$ and $G$ in $L^{p(\cdot)}(w)$ such that
\begin{equation}\label{aF}
\|F\|_{L^{\p}(w)}\leq2,
\end{equation}
\begin{equation}\label{bF}
\|f\|_{L^{\p}(w)}=\|fF^{-1}\|_{L^\infty},
\end{equation}
\begin{equation}\label{aG}
\|G\|_{L^{\p}(w)}\leq2,
\end{equation}
and
\begin{equation}\label{bG}
\|g\|_{L^{\p}(w)}=\|gG^{-1}\|_{L^\infty}.
\end{equation}

Since $w$ belongs to $A_{\p}^{\rho}$, by Theorem~\ref{teo: acot M}  there exists $\te>0$ such that $M_\rho^{\theta}$ is bounded in $L^{\p}(w)$. Thus, applying to $h=F+G$ Lemma~\ref{algoritmo} with $r=p$ and $S=M_\rho^{\theta}$ it follows there exists $H\in L^{\p}(w)$, such that
\begin{equation}\label{hc}
h\leq H,
\end{equation}
and
\begin{equation}\label{hd}
\|H\|_{L^{\p}(w)}\leq 2\|h\|_{L^{\p}(w)}.
\end{equation}
Furthermore, from~\eqref{est rho}, it follows that for $\sigma=\theta(N_\rho+1)$ and $z\in B(x_0,r)\subset\R^d$,
\begin{equation*}
\begin{split}
\bigg(1+\frac r{\rho(x_0)}\bigg)^{-\sigma} \frac{1}{|B(x_0,r)|} &\int_{B(x_0,r)}|H(y)|\,dy\\
&\le2^\theta c_\rho^\theta\bigg(1+\frac {2r}{\rho(z)}\bigg)^{-\theta} \frac{2^d}{|B(z,2r)|} \int_{B(z,2r)}|H(y)|\,dy\\
&\le c_\rho^\theta2^{d+\theta}\, M_\rho^\theta H(z).
\end{split}
\end{equation*}

Then, since $M_\rho^\theta H\le CH$ a.e., we have that $H\in A_1^\rho$. Therefore, from \eqref{hd}, \eqref{aF} and \eqref{aG}, and using the hypothesis with $w=H^{-1}$, we obtain that
\begin{equation*}\label{gfinf}
\begin{split}
\|f\|_{L^{\p}(w)}&=\|f H^{-1} H\|_{L^{\p}(w)} \leq\ \|f H^{-1}\|_{L^\infty}\|H\|_{L^{\p}(w)}\leq C\|gH^{-1}\|_{L^\infty},
\end{split}
\end{equation*}
whenever $\|fH^{-1}\|_{L^\infty}<\infty$. In fact by \eqref{hc} and \eqref{bF}, we have
\begin{equation*}
\|fH^{-1}\|_{L^\infty}\leq\|fF^{-1}\|_{L^\infty}=\|f\|_{L^{\p}(w)}<\infty.
\end{equation*}
	
Finally,
\begin{equation*}
\begin{split}
\|f\|_{L^{\p}(w)}&\le C\|gH^{-1}\|_{L^\infty}\le C\|gG^{-1}\|_{L^\infty}= C\|g\|_{L^{\p}(w)},
\end{split}
\end{equation*}
where we have used \eqref{hc} and \eqref{bG}.

\end{proof}

\begin{cor}\label{cor:ext inf r}
Let  $r>1$ and $\fa=\{(f,g)\}$  a family of pairs of measurable and non-negative functions. Suppose
\begin{equation*}
\|f w\|_{L^\infty} \le C \|g w\|_{L^\infty},
\end{equation*}
holds for every weight $w$ such that $w^{-r} \in A_1^\rho$ and every pairs $(f,g)\in\mathscr{F}$. Then,
\begin{equation}\label{des infpr}
\|f\|_{L^{\p}(w)}\le C \|g\|_{{L^{\p}(w)}},
\end{equation}
holds for every pair $(f,g)\in\fa$, every $p\in\mathcal{P}^{\log}(\R^d)$ with $p^{-}>r$ and every weight $w$ such that $w^r\in A_{\p/r}^{\rho}$, provided that the left hand side of \eqref{des infpr} is finite.
\end{cor}
\begin{proof}
Let  $w$ be such that $w^{-1}\in A_1^\rho$ and $(f,g)\in\mathscr{F}$. Then, using the hypothesis with $w^{1/r}$ it follows
\begin{equation*}
\|f^rw\|^{1/r}_{L^{\infty}}=\|fw^{1/r}\|_{L^{\infty}}\le C\|gw^{1/r}\|_{L^{\infty}}=C\|g^rw\|^{1/r}_{L^{\infty}}.
\end{equation*}

Therefore, we have proved~\eqref{des_inf} for the family $\widetilde{\mathscr{F}}$ of pairs $(f^r,g^r)$ with $(f,g)\in\mathscr{F}$ and any $w$ be such that $w^{-1}\in A_1^\rho$. This is the hypothesis of Theorem~\ref{thm:ext inf}, and applying this theorem, we get
\begin{equation}\label{des infqr}
\|f^r\|_{L^{q(\cdot)}(w)}\le C \|g^r\|_{{L^{q(\cdot)}(w)}},\ \ \ \ (f,g)\in\mathscr{F},
\end{equation}
for every $q\in\mathcal{P}^{\log}(\R^d)$ with $q^{-}>1$ and every $w\in A^\rho_{q(\cdot)}$.

Finally, given $p\in\mathcal{P}^{\log}(\R^d)$ with $p^{-}>r$, we get that $q(\cdot)=\p/r$ belongs to $\mathcal{P}^{\log}(\R^d)$ and $q^{-}>1$. Then, it follows from~\eqref{des infqr} and Lemma~\ref{lem:dilatacion} that for any $w$ that such $w^r\in A^\rho_{q(\cdot)}$, holds that
\begin{equation*}
\begin{split}
\|f\|_{L^{\p}(w)}=\|fw\|_{L^{\p}}&=\|f^rw^r\|^{1/r}_{\frac{\p}r}
=\|f^r\|^{1/r}_{L^{q(\cdot)}(w^r)}\\
&\le C\|g^r\|^{1/r}_{L^{q(\cdot)}(w^r)}=\|g^rw^r\|^{1/r}_{\frac{\p}r}
=\|gw\|_{L^{\p}}=\|g\|_{L^{\p}(w)}.
\end{split}
\end{equation*}
\end{proof}

\section{Bounded mean oscillation results}	

Given a critical radius function $\rho$ and weight $w$  we say that a locally integrable function $f$ belongs to $BMO_{\rho}(w)$ if it satisfies
\begin{equation}\label{osc acot}
\frac{\|\chi_Bw\|_{L^\infty}}{|B|}\int_B|f-f_B|\leq C,\hspace{0.5cm}\text{for all}\ \ B\in\mathcal{B}_\rho,
\end{equation}
and
\begin{equation}\label{prom acot}
\frac{\|\chi_{B(x,\rho(x))}w\|_{L^\infty}}{|{B(x,\rho(x))}|}\int_{B(x,\rho(x))}|f|\leq C,\hspace{0.5cm}\text{for all}\ \ x\in\R^d.
\end{equation}
Here, as usual, $f_B$ stands for the average of $f$ over the ball $B$. A seminorm $\|f\|_{BMO_{\rho}(w)}$ in $BMO_{\rho}(w)$ is defined as the least constant satisfying \eqref{osc acot} and \eqref{prom acot}.

If we take the limit case $\rho\equiv\infty$, the above definition agrees with that of the weighted bounded mean oscillation space considered by Harboure, Mac\'ias y Segovia in~\cite{Pola}.

Given $f\in L^1_\loc(\R^n)$, we define a localized version of the sharp maximal function as

\begin{equation*}
M_\rho^{\sharp} f(x)=\sup_{x\in B\in\mathcal{B}_{\rho}}\frac1{|B|}\int_B|f(y)-f_B|dy\ +\sup_{x\in B=B(z,\rho(z))}\frac1{|B|}\int_B|f(y)|dy.
\end{equation*}

As it is expected, the space $BMO_\rho(w)$ can be described by means of $M_\rho^{\sharp}$ defined above.

\begin{lem}[\cite{BCH-Extrapol}, Lemma 2]\label{lema sharp-bmo}
Let $w$ be a weight. If $f$ belongs to $L^1_\loc(\R^d)$, then
\begin{equation*}
\|f\|_{BMO_{\rho}(w)}\simeq\|w\,M_\rho^{\sharp}f\|_{L^\infty}.
\end{equation*}
\end{lem}

The following proposition gives a $\rho$-local version of Lerner's inequality (see Theorem 1 in \cite{Lerner}).

\begin{prop}[\cite{BCH-Lerner}, Theorem 1]\label{prop: Lerner}
If $f$ and  $u$ are non-negative functions belonging to $L^1_\loc(\R^d)$, then
\begin{equation*}
\int_{\R^d}f(x)u(x)dx\leq C\int_{\R^d}M_{\rho}^\sharp f(x)M_{\rho}u(x)dx,
\end{equation*}
where the constant $C$ is independent of $f$ and $u$.
\end{prop}	

In what follows we present a Fefferman--Stein type inequality relating the local sharp maximal  $M_\rho^{\sharp}$ with $M_\rho$. The precise statement is as follows.

\begin{thm}\label{teo:F-S}
Let $p\in\mathcal{P}^{\log}(\R^d)$  with $p^{-}>1$ and  $w\in A_{\p}^{\rho,\loc}$, then there exist constants $C>0$ and $\beta>1$ such that
\begin{equation}
\|w\,M_\rho f\|_{\p}\le C\|w\,M^\sharp_{\beta\rho}f\|_{\p}
\end{equation}	
\end{thm}

\

To prove this theorem, we need the following point estimate for the sharp maximal function of $M_\rho$.

\begin{prop}[\cite{BCH-Lerner}, Proposition 2]\label{propo:des punt}
Let $0<\delta<1$. There exist constants $C>0$ and $\be>1$ such that
\begin{equation*}
[M_\rho^{\sharp}(M_\rho f)^\delta(x)]^{1/\delta}\le CM_{\beta\rho}^{\sharp}f(x),
\end{equation*}
for every $f\in L^1_\loc(\R^d)$, with $\be$ depending only on the constants $c_\rho$ and $N_\rho$ in~\eqref{est rho}.
\end{prop}

\

\begin{proof}[Proof of Theorem~\ref{teo:F-S}]
Let $0<\delta<1$ and $r(\cdot)=\p/\delta$, then we get that $r$ and $r'$ belongs to $\mathcal{P}^{\log}(\R^d)$. By Lemma~\ref{lem: norma dual} we have 
\begin{equation*}
\|w\,M_\rho f\|_{\p}=\|w^\delta (M_\rho f)^\delta\|^{1/\delta}_{r(\cdot)}
\lesssim\sup_{\|w^{-\delta}g\|_{r'(\cdot)}\le1}\left(\int_{\R^d}(M_\rho f)^\delta|g|\,dx\right)^{1/\delta}.
\end{equation*}

Then by applying Proposition~\ref{prop: Lerner} and the generalized H\"older inequality  we obtain that
\begin{equation*}
\begin{split}
\|w\,M_\rho f\|_{\p}&\lesssim\sup_{\|w^{-\delta}g\|_{r'(\cdot)}\le1}\left(\int_{\R^d}M_\rho^{\sharp}((M_\rho f)^\delta)M_\rho g\,dx\right)^{1/\delta}\\
&\le\sup_{\|w^{-\delta}g\|_{r'(\cdot)}\le1}\|w^\delta M_\rho^{\sharp}(M_\rho f)^\delta\|^{1/\delta}_{r(\cdot)}\ \|w^{-\delta}M_\rho g\|_{r'(\cdot)}^{1/\delta}\\
&=\|w [M_\rho^{\sharp}(M_\rho f)^\delta]^{1/\delta}\|_{\p}\ \sup_{\|w^{-\delta}g\|_{r'(\cdot)}\le1}\|w^{-\delta}M_\rho g\|_{r'(\cdot)}^{1/\delta}.
\end{split}
\end{equation*}

Since $w\in A_{\p}^{\rho,\loc}$ it follows that $w^\delta\in A_{r(\cdot)}^{\rho,\loc}$ (see Proposition~\ref{prop: wdelta}), which is equivalent to $w^{-\delta}\in A_{r'(\cdot)}^{\rho,\loc}$. Then, by Theorem~\ref{teo: acot M} and Proposition~\ref{propo:des punt} we obtain that
\begin{equation*}
\begin{split}
\|w\,M_\rho f\|_{\p}&\lesssim\|w\,M_{\beta\rho}^{\sharp}f\|_{\p}\ \sup_{\|w^{-\delta}g\|_{r'(\cdot)}\le1}\|w^{-\delta} g\|_{r'(\cdot)}^{1/\delta}\le\|w\,M_{\beta\rho}^{\sharp}f\|_{\p}.
\end{split}
\end{equation*}
\end{proof}

\begin{cor}\label{cor:F-S}
Let $p\in\mathcal{P}^{\log}(\R^d)$  with $p^{-}>1$ and $w\in A_{\p}^{\rho,\loc}$. Then there is $\gamma<1$ such that
\begin{equation*}
\|w\,M_{\gamma\rho} f\|_{\p}\le C\|w\,M^\sharp_{\rho}f\|_{\p}
\end{equation*}
\end{cor}
\begin{proof}
Let $\be>1$. Observe that for $\gamma=\frac1\beta<1$, the critical radius functions $\rho$ and $\gamma\rho$ satisfy inequality~\eqref{est rho} with the same constants ($c_{\rho}$ and $N_\rho$). Also, according to Proposition~\ref{prop: Aploc gamma}, classes $A_{\p}^{\rho,\loc}$ and $A_{\p}^{\gamma\rho,\loc}$ are the same. Therefore, Theorem~\ref{teo:F-S} is true exchanging $\rho$ by $\gamma\rho$.

\end{proof}

As a consequence of Corollary~\ref{cor:F-S} and the fact that almost everywhere $g\le M_{\gamma\rho}g$ for every $g\in L^1_{\loc}(\R^d)$ and $\gamma>0$, we have the following corollary.

\begin{cor}\label{des norma f y f sharp}
If $p\in\mathcal{P}^{\log}(\R^d)$  with $p^{-}>1$, $w\in A_{\p}^{\rho,\loc}$ and $g$ belongs to $L^1_{\loc}(\R^d)$, there exists a constant $C>0$ such that
\begin{equation*}
\|g\,w\|_{\p}\leq C\|M_\rho^{\sharp} g\,w\|_{\p}.
\end{equation*}
\end{cor}

From these results we can reformulate our extrapolation result to include the case of operators that are bounded from $L^\infty(w)$ into $BMO_\rho(w)$.

\begin{thm}\label{thm:ext_bmoT}
Let $r\ge1$ and $T$ a bounded operator from $L^\infty(w)$ into $BMO_\rho(w)$ for any $w$ such that $w^{-r}\in A^\rho_1$. Then
\begin{equation*}
\|T f\|_{L^{\p}(w)}\le C\|f\|_{L^{\p}(w)},
\end{equation*}
holds for  every $p\in\mathcal{P}^{\log}(\R^d)$ with $p^{-}>r$ and every weight $w$ such that $w^r\in A_{\p/r}^{\rho}$.
\end{thm}
\begin{proof}
By the hypothesis on $T$ and Lemma~\ref{lema sharp-bmo} we have
\begin{equation*}
\|(M^\sharp_{\rho}\circ T)f\,w\|_{L^\infty}\leq C\|f\,w\|_{L^\infty},
\end{equation*}
for every weight $w$ such that $w^{-r}\in A^\rho_1$. We can then apply Theorem~\ref{thm:ext inf} or Corollary~\ref{cor:ext inf r} to conclude that $M^\sharp_{\rho}\circ T$ is bounded on $L^{\p}(w)$ for every weight $w$ such that $w^r\in A_{\p/r}^{\rho}$. Hence, the boundedness of $T$ follows from Corollary~\ref{des norma f y f sharp} since by Proposition~\ref{prop: wdelta} $w\in  A_{\p}^{\rho}\subset  A_{\p}^{\rho,\loc}$.
	
\end{proof}

\section{Schr\"odinger Type Singular Integrals}

In~\cite{BCH-Extrapol} a class of operators resembling those of Calder\'on--Zygmund theory, but adapted to the Schr\"odinger context, was introduced (see also~\cite{BCH-Lerner} and~\cite{BCH-Singular Int}).  In it the authors analyzed their behavior on weighted $L^p$ spaces for $1<p\le\infty$ using extrapolation techniques. This type of operators were also considered in~\cite{BCQ-Weighted STS} and~\cite{Ma-Stinga-et}, where conditions are given to obtain their boundedness in regularity spaces in the context of Schr\"odinger with and without weights, respectively.

We shall call  Schr\"odinger--Calder\'on--Zygmund operator of type $(\infty,\delta)$ for $0<\delta\le1$ to an operator $T$ such that
\begin{enumerate}[(I)]
\item\label{cond I} $T$ is bounded in $L^p$ for some $1<p<\infty$.
\item $T$ has an associated kernel $K:\R^d\times\R^d\longrightarrow\R$, in the sense that
\begin{equation*}
Tf(x)=\int_{\R^d}K(x,y)f(y)\,dy,\ \ \ \ f\in L_c^{\infty}\ \ \ \text{and\ \  a.e.}\  x\notin\text{supp}\, f.
\end{equation*}
Further, for each $N>0$ there exists a constant $C_N>0$ such that
\begin{equation}\label{cond tamano nucleo T para s infinito}
|K(x,y)|\leq C_N\frac1{|x-y|^d}\bigg(1+\frac{|x-y|}{\rho(x)}\bigg)^{-N},
\end{equation}
for any  $x\neq y$, and there exists $C>0$ such that
 \begin{equation}\label{cond suavidad nucleo T para s infinito}
|K(x,y)-K(x_0,y)|\leq C\frac{|x-x_0|^\delta}{|x-y|^{d+\delta}},
\end{equation}
for every $x,y\in \R^d$, whenever $|x-x_0|<\frac{|x-y|}2$.
\end{enumerate}

\

Operators of this kind, such the Schr\"odinger Riesz transforms $\nabla\mathcal{L}^{-1/2}$, appear in this context when the potential $V$ satisfies a reverse--H\"older condition~\eqref{reverse} for some $q\ge d$. Nevertheless, when $d/2<q<d$,  due to the lack of regularity of the potential $V$, some important operators related to Schr\"odinger semi-group do not share the above properties.

To deal with those cases, following~\cite{BCQ-Weighted STS}, we introduce the following class. We shall say that a linear operator
$T$ is a Schr\"odinger--Calder\'on--Zygmund operator of type $(s,\delta)$, for $1<s<\infty$ and $0<\delta\le1$,  if
\begin{enumerate}[(I$_{s}$)]
\item $T$ is bounded from $L^{s'}$ into $L^{s',\infty}$.
\item $T$ has an associated kernel $K:\R^d\times\R^d\longrightarrow\R$, in the sense that
\begin{equation*}
Tf(x)=\int_{\R^d}K(x,y)f(y)\,dy,\ \ \ \ f\in L_c^{\infty}\ \ \text{and\ \ a.e.}\  x\notin\text{supp}\, f.
\end{equation*}
Further, for each $N>0$ there exists a constant $C_N>0$ such that
\begin{equation}\label{cond tamano nucleo T para s}
\left(\frac1{R^d}\int_{R<|x_0-y|<2R}|K(x,y)|^{s}dy\right)^{1/s}\leq C_NR^{-d}\bigg(1+\frac{R}{\rho(x)}\bigg)^{-N},
\end{equation}
for any  $|x-x_0|<R/2$, and there exists $C>0$ such that
\begin{equation}\label{cond suavidad nucleo T para s}
\left(\frac1{R^d}\int_{R<|x_0-y|<2R}|K(x,y)-K(x_0,y)|^{s}dy\right)^{1/s}\leq CR^{-d}\left(\frac{r}{R}\right)^\delta,
\end{equation}
for every $|x-x_0|<r<\rho(x_0)$ and $r<R/2$.
\end{enumerate}

\begin{rem}
Is straightforward to check that if a kernel satisfies~\eqref{cond tamano nucleo T para s} and~\eqref{cond suavidad nucleo T para s}  for some $s>1$, then it also satisfies both conditions for any $1\le t<s$.
\end{rem}

\begin{rem}\label{imp}
It is clear that if $T$ is a Schr\"odinger--Calder\'on--Zygmund operator of type $(\infty,\delta)$ which verifies condition~\eqref{cond I} with $p=p_0$, then $T$ is a Schr\"odinger--Calder\'on--Zygmund operator of type $(p_0,\delta)$, since the  conditions ~\eqref{cond tamano nucleo T para s} and~\eqref{cond suavidad nucleo T para s}  follow easily from the pointwise estimates of the kernel required in~\eqref{cond tamano nucleo T para s infinito} and~\eqref{cond suavidad nucleo T para s infinito}. 
\end{rem}

From the results of~\cite{BCH-Extrapol} and~\cite{BCH-Singular Int} we have the following theorem concerning strong weighted inequalities for Schr\"odinger--Calder\'on--Zygmund operator of type $(\infty,\delta)$ and $(s,\delta)$.

\begin{thm}
Let $0<\delta\le1$.
\begin{enumerate}
\item If $T$ is a Schr\"odinger--Calder\'on--Zygmund operator of type $(\infty,\delta)$,  then $T$ is bounded on $L^p(w)$ for every $1<p<\infty$ and any $w\in A_{p}^{\rho}$.
\item If $T$ is a Schr\"odinger--Calder\'on--Zygmund operator of type $(s,\delta)$,  then $T$ is bounded on $L^p(w)$ for every $s'<p<\infty$ and any $w\in A_{p/s'}^{\rho}$.
\end{enumerate}
\end{thm}

\begin{proof}
The second statement is a direct consequence of Proposition 6 and Theorem 5 of~\cite{BCH-Extrapol}.
	
Regarding the first statement, as indicated in Remark~\ref{imp}, if $T$ is a Schr\"odinger--Calder\'on--Zygmund operator of type $(\infty,\delta)$ which verifies condition~\eqref{cond I} with $p=p_0$, then $T$ is a Schr\"odinger--Calder\'on--Zygmund operator of type $(p_0,\delta)$. Moreover, as we have mentioned, its follows from  Proposition 6 and Theorem 5 in~\cite{BCH-Extrapol} 
that  such operators are bounded on $L^p(w)$ for all $p'_0<p<\infty$ and every
$w\in A^\rho_{p/p'_0}$. This implies, in particular, that such operators are bounded in $L^p$ for all $p'_0<p<\infty$.

Under these conditions, i.e. an operator $T$ that is bounded on $L^p$ for any $p>q$, for some $q>1$ and whose kernel verifies~\eqref{cond tamano nucleo T para s infinito} and~\eqref{cond suavidad nucleo T para s infinito},  was proved in Proposition 3.2 in~\cite{BCH-Singular Int} that $T$ is bounded on $L^p(w)$ for all $1<p<\infty$ and any $w\in A^\rho_{p}$.
 
\end{proof}

\begin{rem}\label{cond equiv}
From the above it follows that a Schr\"odinger--Calder\'on--Zygmund operator of type $(\infty,\delta)$ is bounded in $L^p$ for all $1<p<\infty$. Therefore condition~\eqref{cond I} in the class definition $(\infty,\delta)$ is
equivalent to ask boundedness on $L^p$ for any $1<p<\infty$.
	
We also get that if $T$ is a Schr\"odinger--Calder\'on--Zygmund operator of type $(\infty,\delta)$, then $T$ is a Schr\"odinger--Calder\'on--Zygmund operator of type $(s,\delta)$, for any $1<s<\infty$.
\end{rem}

As we have already mentioned, our main aim is to prove boundedness results of such operators on weighted variable Lebesgue spaces  $L^{\p}(w)$ with weights defined in terms of $\rho$. We start with the following result.

\begin{thm}\label{thm: SCZ inf}
Let  $p\in\mathcal{P}^{\log}(\R^d)$ with $p^{-}>1$. If $T$ is a Schr\"odinger--Calder\'on--Zygmund operator of type $(\infty,\delta)$, then
\begin{equation*}
\|T f\|_{L^{\p}(w)}\le C\|f\|_{L^{\p}(w)},
\end{equation*}
holds for  every weight  $w\in A_{\p}^{\rho}$. Moreover, its adjoint operator $T^{*}$ is also bounded on $L^{\p}(w)$  for every $w\in A_{\p}^{\rho}$.
\end{thm}

An important element in the proof of the above theorem is given by the following result.

\begin{prop}
If $T$ is a Schr\"odinger--Calder\'on--Zygmund operator of type $(\infty,\delta)$, $0<\delta\le1$, then $T$ is bounded from $L^\infty(w)$ into $BMO_\rho(w)$ for every weight $w$ such that $w^{-1}\in A_1^\rho$.
\end{prop}
\begin{proof}
By the Proposition 5 in~\cite{BCH-Extrapol} a bounded operator on $L^p$ for any $1<p<\infty$ whose kernel verifies conditions~\eqref{cond tamano nucleo T para s infinito} and~\eqref{cond suavidad nucleo T para s infinito} is bounded from $L^\infty(w)$ into $BMO_\rho(w)$ for every weight $w$ such that $w^{-1}\in A_1^\rho$. According to this and the Remark~\ref{cond equiv}, the thesis is followed to for any Schr\"odinger--Calder\'on--Zygmund operator of type $(\infty,\delta)$.
	
\end{proof}

\begin{proof}[Proof of Theorem~\ref{thm: SCZ inf}]
The first statement follows directly from Theorem~\ref{thm:ext_bmoT} and the above proposition.

Regarding the second statement, by Lemma~\ref{lem: norma dual} and Lemma~\ref{lem: Holder} we obtain that
\begin{equation*}
\begin{split}
\|T^*f\|_{L^{\p}(w)}=\|w\,T^*f\|_{\p}
&\lesssim\sup_{\|w^{-1}g\|_{\pri}\le1}\ \left|\int_{\R^d}T^*f(x)\,g(x)\,dx\right|\\
&=\sup_{\|w^{-1}g\|_{\pri}\le1}\ \left|\int_{\R^d}f(x)\,Tg(x)\,dx\right|\\
&\lesssim\sup_{\|w^{-1}g\|_{\pri}\le1}\ \|w\,f\|_{\p}\,\|w^{-1}Tg\|_{\pri}\\
&\le\|w\,f\|_{\p}\,\sup_{\|w^{-1}g\|_{\pri}\le1}\,\|w^{-1}Tg\|_{\pri}.
\end{split}
\end{equation*}

Since that for $p\in\mathcal{P}^{\log}(\R^d)$, it is $p^+<\infty$, it follows that $(p')^{-}=(p^+)'>1$. Also, condition $w\in A_{\p}^{\rho}$ implies that $w^{-1}\in A_{\pri}^{\rho}$. According to what has already been proven in the first part,  it can be concluded that
\begin{equation*}
\begin{split}
\|T^*f\|_{L^{\p}(w)}&\lesssim\|w\,f\|_{\p}\,\sup_{\|w^{-1}g\|_{\pri}\le1}\,\|w^{-1}Tg\|_{\pri}\\
&\lesssim\|w\,f\|_{\p}\,\sup_{\|w^{-1}g\|_{\pri}\le1}\,\|w^{-1}g\|_{\pri}\\
&=\|w\,f\|_{\p}.
\end{split}
\end{equation*}

\end{proof}

Let us now see the corresponding result for the operators of type $(s,\delta)$.

\begin{thm}\label{thm: SCZ s}
Let $s>1$ and $p\in\mathcal{P}^{\log}(\R^d)$ with $p^{-}>s'$. If $T$ is a Schr\"odinger--Calder\'on--Zygmund operator of type $(s,\delta)$, then
\begin{equation*}
\|T f\|_{L^{\p}(w)}\le C\|f\|_{L^{\p}(w)},
\end{equation*}
holds for  every weight $w$ such that $w^{s'}\in A_{\p/s'}^{\rho}$. Moreover, its adjoint operator $T^{*}$ is bounded on $L^{\p}(w)$  for every $p\in\mathcal{P}^{\log}(\R^d)$ with $1<p^{+}<s$ and any $w$ such that $w^{-s'}\in A_{\pri/s'}^{\rho}$.
\end{thm}

A key element in the proof of the previous theorem is given by the next proposition.

\begin{prop}
Let $s>1$ and $T$ a Schr\"odinger--Calder\'on--Zygmund operator of type $(s,\delta)$, $0<\delta\le1$, then $T$ is bounded from $L^\infty(w)$ into $BMO_\rho(w)$ for every weight $w$ such that $w^{-s'}\in A_1^\rho$.
\end{prop}
\begin{proof}
According to Proposition 6 in~\cite{BCH-Extrapol}, an operator bounded from $L^{s'}$ into $L^{s',\infty}$ and whose kernel satisfies the following conditions is bounded from $L^\infty(w)$ into $BMO_\rho(w)$ for every weight $w$ such that $w^{-s'}\in A_1^\rho$, provided that

\begin{enumerate}[i)]
\item for each $N>0$ there exists a constant $C_N$ such that
\begin{equation}\label{cond tamano integral}
\left(\int_{R<|x_0-y|<2R}|K(x,y)|^{s}dy\right)^{1/s}\leq C_NR^{-d/s'}\bigg(\frac{\rho(x_0)}{R}\bigg)^{N},
\end{equation}
for every $x\in B(x_0,\rho(x_0))$ and $R>2\rho(x_0)$;
\item there exists a constant $C$ such that
\begin{equation}\label{cond suavidad integral}
\sum_{k\ge1}(2^kr)^{d/s'}\left(\int_{B_{k+1}\setminus B_k}|K(x,y)-K(x_0,y)|^{s}dy\right)^{1/s}\leq C,
\end{equation}
for every $x\in B(x_0,r)$ con $r\le\rho(x_0)$ and $B_k=B(x_0,2^kr)$, $k\in\N$.
\end{enumerate}

Therefore, we need to check that conditions~\eqref{cond tamano nucleo T para s} and~\eqref{cond suavidad nucleo T para s}  implies~\eqref{cond tamano integral} and~\eqref{cond suavidad integral} respectively.

On the one hand, condition~\eqref{cond tamano integral} follows directly from~\eqref{cond tamano nucleo T para s} taking into account that from~\eqref{est rho}, for each $x\in B(x_0,\rho(x_0))$,  we have to $\rho(x)\simeq\rho(x_0)$, more specifically $(c_\rho2^{N_\rho})^{-1}\rho(x_0)\le\rho(x)\le c_\rho2^{\frac{N_\rho}{N_\rho+1}}\rho(x_0)$.

On the other hand, from~\eqref{cond suavidad nucleo T para s}, for every $x\in B(x_0,r)$ con $r\le\rho(x_0)$, we get
\begin{equation*}
\begin{split}
\sum_{k\ge1}(2^kr)^{d/s'}&\left(\int_{B_{k+1}\setminus B_k}|K(x,y)-K(x_0,y)|^{s}dy\right)^{1/s}\\
&\le C\sum_{k\ge1}(2^kr)^{d/s'+d/s-d}2^{-k\delta}\le C.
\end{split}
\end{equation*}
\end{proof}

\begin{proof}[Proof of Theorem~\ref{thm: SCZ s}]
The first statement follows directly from Theorem~\ref{thm:ext_bmoT} and the above proposition.

On the other hand, analogously to the proof of Theorem~\ref{thm: SCZ inf}, we have that
\begin{equation*}
\|T^*f\|_{L^{\p}(w)}=\|w\,T^*f\|_{\p}\lesssim\|w\,f\|_{\p}\,\sup_{\|w^{-1}g\|_{\pri}\le1}\ \|w^{-1}Tg\|_{\pri}
\end{equation*}

Since that for $p\in\mathcal{P}^{\log}(\R^d)$ with $1<p^{+}<s$ it is verified that $(p')^{-}=(p^+)'>s'$, then according to what has already been proven in the first part, for all weight $w$ such that $w^{-s'}\in A_{\pri/s'}^{\rho}$ it follows that
\begin{equation*}
\begin{split}
\|T^*f\|_{L^{\p}(w)}&\lesssim\|w\,f\|_{\p}\,\sup_{\|w^{-1}g\|_{\pri}\le1}\ \|w^{-1}Tg\|_{\pri}\\
&\lesssim\|w\,f\|_{\p}\,\sup_{\|w^{-1}g\|_{\pri}\le1}\ \|w^{-1}g\|_{\pri}\\
&=\|w\,f\|_{\p}.
\end{split}
\end{equation*}

\end{proof}

\section{Application to Schr\"odinger operators}

In this  section we go back to the setting of the Schr\"odinger operator
$\mathcal{L}=-\Delta+V$. 
We will apply Theorems~\ref{thm: SCZ inf} and~\ref{thm: SCZ s} to obtain continuity  on $L^{\p}(w)$ for  for  several  operators related to $\mathcal{L}$ with potentials satisfying a reverse H\"older condition as stated in the introduction. Therefore, from now on, $\rho$ will be the critical radius function associated to the potential by means of~\eqref{rho Schrodinger}.

For shortness we will use the notation SCZ to mean a Schr\"odinger--Calder\'on--Zygmund operator when $\rho$ is the critical radius function derived from the potential $V$.

\subsection{Application to Schr\"odinger--Riesz Transforms}
In this section we consider the singular integral operators known as the Riesz--Schr\"odinger transforms of order 1 and 2 given by $R_1=\nabla\mathcal{L}^{-1/2}$ and $R_2=\nabla^2\mathcal{L}^{-1}$
respectively, together with their adjoints $R^*_1=\mathcal{L}^{-1/2}\nabla$ and $R^*_2=\mathcal{L}^{-1}\nabla^2$.

Let $V\in RH_q$ with $q>d/2$. We will analyze the operators $R_1$, $R^*_1$, $R_2$ and $R^*_2$  considering the cases $q\ge  d$ and $d/2<q<d$.

\begin{thm}
Let $V\in RH_q$ with $q\ge d$. Then,  the operators $R_1$ and $R_1^*$ are SCZ operators of type $(\infty,\delta)$. Moreover, $R_1$ and $R_1^*$ are bounded on $L^{\p}(w)$ for every $p\in\mathcal{P}^{\log}(\R^d)$ with $p^{-}>1$ and any $w\in A_{\p}^{\rho}$.
\end{thm}
\begin{proof}
It was proven in~\cite{Shen} (see Theorem 0.8) that  $R_1$ and $R_1^*$ are Calder\'on--Zygmund operators. So, in that case, $R_1$ and $R_1^*$ are bounded on $L^p$ for all $1<p<\infty$ and their kernels satisfy~\eqref{cond suavidad nucleo T para s infinito}. Moreover, their kernels also satisfy~\eqref{cond tamano nucleo T para s infinito} (see estimate (6.5)  given there).
Therefore, both $R_1$ and $R_1^*$ are Schr\"odinger--Calder\'on--Zygmund operators of
type $(\infty,\delta)$. However, they have different degree of smoothness $\delta=1-d/q$ in the
first case and $\delta=1$ in the second case. Finally, the second statement follows directly from Theorem~\ref{thm: SCZ inf}.

\end{proof}

\begin{thm}
Let $V\in RH_q$ with $d/2<q<d$. Then,  the operator $R_1^*$ is SCZ operator of type $(s,\delta)$ with $s$ such that $1/s=1/q-1/d$.

Moreover,  $R_1^*$ is bounded on $L^{\p}(w)$ for every $p\in\mathcal{P}^{\log}(\R^d)$ with $p^{-}>s'$ and  every weight $w$ such that $w^{s'}\in A_{\p/s'}^{\rho}$, while $R_1$ is bounded on $L^{\p}(w)$ for every $p\in\mathcal{P}^{\log}(\R^d)$ with $1<p^{+}<s$ and  every weight $w$ such that $w^{-s'}\in A_{\pri/s'}^{\rho}$.
\end{thm}
\begin{proof}
Theorem 0.5 in~\cite{Shen} gives us $R_1^*$ is bounded in $L^p$ for $s'\le p<\infty$ with $s$ such that $1/s=1/q-1/d$. The proof that the kernel of  $R_1^*$ satisfy conditions~\eqref{cond tamano nucleo T para s} and~\eqref{cond suavidad nucleo T para s} it is contained in~\cite{BCQ-Weighted STS} (see (77) there). Then, the bounding of the operators  $R_1$ and $R_1^*$ in $L^{\p}(w)$ follows directly from Theorem~\ref{thm: SCZ s}.

\end{proof}

\begin{thm}
Let $V\in RH_q$ with $q>d/2$. Then,  the operator $R^*_2$ is SCZ operator of type $(q,\delta)$ with $\delta=\min\{1,2-d/q\}$.
	
Moreover,  $R_2^*$ is bounded on $L^{\p}(w)$ for every $p\in\mathcal{P}^{\log}(\R^d)$ with $p^{-}>q'$ and  every weight $w$ such that $w^{q'}\in A_{\p/q'}^{\rho}$, while $R_2$ is bounded on $L^{\p}(w)$ for every $p\in\mathcal{P}^{\log}(\R^d)$ with $1<p^{+}<q$ and  every weight $w$ such that $w^{-q'}\in A_{\pri/q'}^{\rho}$.
\end{thm}
\begin{proof}
Firstly, it follows from Theorem 0.3 in~\cite{Shen}, that $R_2^*$ is bounded on $L^{p}$ for $q'\le p<\infty$. The proof that the kernel of  $R_2^*$ satisfy conditions~\eqref{cond tamano nucleo T para s} and~\eqref{cond suavidad nucleo T para s} is contained in the Proposition 8 of~\cite{BCQ-Weighted STS}. Therefore, applying Theorem~\ref{thm: SCZ s} we complete the proof.

\end{proof}

In~\cite{BCQ-Behaviour of Schrodinger} the authors show that if $V\in RH_q$ with $q>d/2$ and verifies  a local smoothness condition relative to the function $\rho$, namely 
\begin{equation}\label{suav local V}
|V(x)-V(y)|\lesssim\frac{|x-y|^\alpha}{\rho(x)^{\alpha+2}},
\end{equation}
for every $x,y\in\R$ such that $|x-y|<\rho(x)$  and some $0<\alpha\le1$, then $R_2$ is a Schr\"odinger--Calder\'on--Zygmund operator of type $(\infty,\alpha)$ (see Proposition 3 there).

As a consequence of this and Theorem~\ref{thm: SCZ inf} we have the following result.
\begin{thm}
Let $V\in RH_q$ with $q>d/2$ and suppose that $V$ satisfies~\eqref{suav local V} for some $0<\alpha\le1$ . Then,  the operators  $R_2$ and $R_2^*$ are bounded on $L^{\p}(w)$ for every $p\in\mathcal{P}^{\log}(\R^d)$ with $p^{-}>1$ and any $w\in A_{\p}^{\rho}$.
\end{thm}

\subsection{Applications to Schr\"odinger--Riesz transforms involving $V$}
In this section we consider the operators $M_\gamma=\mathcal{L}^{-\gamma}V^{\gamma}$ for $0<\gamma<d/2$  and $N_\gamma=\mathcal{L}^{-\gamma}\nabla V^{\gamma-1/2}$ for $1/2<\gamma\le1$. These operators have been first considered by Shen in~\cite{Shen}, the first one for the cases $\gamma=1/2$ and $\gamma=1$ and the second one only for the case $\gamma=1$.

In~\cite{BCQ-Weighted STS} the authors have proved the following results directly related to the operators in question.

\begin{prop}[\cite{BCQ-Weighted STS}, Propositions 6 and 7]
Let $V\in RH_q$ with $q>d/2$. Then,
\begin{enumerate}
\item The operator $N_\gamma$  is a SCZ  operator of type $(s,\delta)$  for $1/2<\gamma\le1$, with $\delta=\{1,2-d/q\}$ and $s$ such that $\frac1s=\left(\frac1q-\frac1d\right)^{+}+\frac{2\gamma-1}{2q}$, where $\left(\frac1q-\frac1d\right)^{+}=\max\big\{\frac1q-\frac1d,0\big\}$.
\item The operator $M_\gamma$  is a SCZ  operator of type $(s,\delta)$  for $0<\gamma<d/2$, with $\delta<\{1,2-d/q\}$  and $s=q/\gamma$.
\end{enumerate}
\end{prop}

As an application of Theorem~\ref{thm: SCZ s} and the above proposition we get boundedness properties for these operators on $L^{\p}(w)$.

\begin{thm}
Let $V\in RH_q$ with $q>d/2$, $1/2<\gamma\le1$ and $s$ such that $\frac1s=\left(\frac1q-\frac1d\right)^{+}+\frac{2\gamma-1}{2q}$. Then,  the operator $N_\gamma$  is  bounded on $L^{\p}(w)$ for every $p\in\mathcal{P}^{\log}(\R^d)$ with $p^{-}>s'$ and  every weight $w$ such that $w^{s'}\in A_{\p/s'}^{\rho}$. Moreover, its adjoint operator $V^{\gamma-1/2}\nabla\mathcal{L}^{-\gamma} $ is bounded on $L^{\p}(w)$  for every $p\in\mathcal{P}^{\log}(\R^d)$ with $1<p^{+}<s$ and any $w$ such that $w^{-s'}\in A_{\pri/s'}^{\rho}$.
\end{thm}

\begin{thm}
Let  $V\in RH_q$ with $q>d/2$, $0<\gamma<d/2$ and $s=q/\gamma$. Then,  the operator $M_\gamma$  is  bounded on $L^{\p}(w)$ for every $p\in\mathcal{P}^{\log}(\R^d)$ with $p^{-}>s'$ and  every weight $w$ such that $w^{s'}\in A_{\p/s'}^{\rho}$. Moreover, its adjoint operator $V^{\gamma}\mathcal{L}^{-\gamma} $ is bounded on $L^{\p}(w)$  for every $p\in\mathcal{P}^{\log}(\R^d)$ with $1<p^{+}<s$ and any $w$ such that $w^{-s'}\in A_{\pri/s'}^{\rho}$.
\end{thm}

Applying these results to the operators considered by Shen, the following results are obtained.

\begin{thm}
Let $V\in RH_q$ with $q>d/2$ and $s$ such that $s=2q$ if $q\ge d$ or $\frac1s=\frac{3}{2q}-\frac1d$ if $d/2<q<d$. Then,  the operator $\mathcal{L}^{-1}\nabla V^{1/2}$  is  bounded on $L^{\p}(w)$ for every $p\in\mathcal{P}^{\log}(\R^d)$ with $p^{-}>s'$ and  every weight $w$ such that $w^{s'}\in A_{\p/s'}^{\rho}$.

Moreover, the operator $V^{1/2}\nabla\mathcal{L}^{-1} $ is bounded on $L^{\p}(w)$  for every $p\in\mathcal{P}^{\log}(\R^d)$ with $1<p^{+}<s$ and any $w$ such that $w^{-s'}\in A_{\pri/s'}^{\rho}$.
\end{thm}

\begin{thm}
Let $V\in RH_q$ with $q>d/2$. Then,  the operator $\mathcal{L}^{-1/2}V^{1/2}$  is  bounded on $L^{\p}(w)$ for every $p\in\mathcal{P}^{\log}(\R^d)$ with $p^{-}>(2q)'$ and  every weight $w$ such that $w^{(2q)'}\in A_{\p/(2q)'}^{\rho}$.
	
Moreover, the operator $V^{1/2}\mathcal{L}^{-1/2} $ is bounded on $L^{\p}(w)$  for every $p\in\mathcal{P}^{\log}(\R^d)$ with $1<p^{+}<2q$ and any $w$ such that $w^{-(2q)'}\in A_{\pri/(2q)'}^{\rho}$.
\end{thm}

\begin{thm}
Let $V\in RH_q$ with $q>d/2$. Then,  the operator $\mathcal{L}^{-1}V$  is  bounded on $L^{\p}(w)$ for every $p\in\mathcal{P}^{\log}(\R^d)$ with $p^{-}>q'$ and  every weight $w$ such that $w^{q'}\in A_{\p/q'}^{\rho}$.
	
Moreover, the operator $V\mathcal{L}^{-1} $ is bounded on $L^{\p}(w)$  for every $p\in\mathcal{P}^{\log}(\R^d)$ with $1<p^{+}<q$ and any $w$ such that $w^{-q'}\in A_{\pri/q'}^{\rho}$.
\end{thm}

\subsection{Applications to operators of the type $\mathcal{L}^{i\alpha}$}
In this section we consider the power operators $\mathcal{L}^{i\alpha}$ for $\alpha\in\R$. From what has been done in~\cite{Shen} it is easy to verify that the following theorem holds.

\begin{thm}
Let $V\in RH_q$ with $q>d/2$. Then,  the operator $\mathcal{L}^{i\alpha}$ with $\alpha\in\R$ is a SCZ operator of type $(\infty,\delta)$. Moreover, $\mathcal{L}^{i\alpha}$ is bounded on $L^{\p}(w)$ for every $p\in\mathcal{P}^{\log}(\R^d)$ with $p^{-}>1$ and any $w\in A_{\p}^{\rho}$.
\end{thm}
\begin{proof}
It was proved in~\cite{Shen} (see Theorem 0.4) that $\mathcal{L}^{i\alpha}$ with $\alpha\in\R$ is a Calder\'on--Zygmund operator. So, in that case, $\mathcal{L}^{i\alpha}$ is bounded on $L^p$ for all $1<p<\infty$ and its kernel satisfy~\eqref{cond suavidad nucleo T para s infinito}. In addition it satisfies the regularity condition~\eqref{cond tamano nucleo T para s infinito} (see equation (4.3) given there). Therefore, applying Theorem~\ref{thm: SCZ inf} we complete the proof.

\end{proof}

\section{Declarations}

\subsection{Ethical Approval} Not applicable.

\subsection{Competing interests} There are no.

\subsection{Authors' contributions } Not applicable.

\subsection{Funding} This research is supported by  Universidad Nacional del Nordeste (UNNE) and Consejo Nacional de Investigaciones Cient\'ificas y T\'ecnicas (CONICET), Argentina.

\subsection{Availability of data and materials } Not applicable.



\end{document}